\documentclass[sn-mathphys, Numbered,onefignum,onetabnum]{sn-jnl}

\usepackage{graphicx}%
\usepackage{tikz}
\usepackage{multirow}%
\usepackage{amsmath,amssymb,amsfonts}%
\usepackage{amsthm}%
\usepackage{mathrsfs}%
\usepackage[title]{appendix}%
\usepackage{xcolor}%
\usepackage{textcomp}%
\usepackage{manyfoot}%
\usepackage{booktabs}%
\usepackage{algorithm}%
\usepackage{algorithmicx}%
\usepackage{algpseudocode}%
\usepackage{listings}%

\usetikzlibrary{arrows,positioning}
\usetikzlibrary{shapes}
\usetikzlibrary{decorations.pathmorphing}
\usetikzlibrary{decorations.pathreplacing}
\usetikzlibrary{decorations.shapes}
\usetikzlibrary{decorations.text}
\usetikzlibrary{decorations.markings}
\usetikzlibrary{snakes}
\usetikzlibrary{patterns}

\newcommand{\R}{\mathbb{R}}
\newcommand{\norm}[1]{\left\lVert#1\right\rVert}
 
\usepackage{lipsum}
\usepackage{epstopdf}
\ifpdf
  \DeclareGraphicsExtensions{.eps,.pdf,.png,.jpg}
\else
  \DeclareGraphicsExtensions{.eps}
\fi

\DeclareMathOperator*{\argmax}{arg\,max}

\usepackage{enumitem}
\setlist[enumerate]{leftmargin=.5in}
\setlist[itemize]{leftmargin=.5in}

\usepackage[normalem]{ulem}

\usepackage{amsopn}

\makeatletter
\newcommand*{\addFileDependency}[1]{
  \typeout{(#1)}
  \@addtofilelist{#1}
  \IfFileExists{#1}{}{\typeout{No file #1.}}
}
\makeatother

\begin{document}

\title[A Preconditioned Interior Point Method for Support Vector Machines Using an ANOVA-Decomposition and NFFT-Based Matrix--Vector Products]{A Preconditioned Interior Point Method for Support Vector Machines Using an ANOVA-Decomposition and NFFT-Based Matrix--Vector Products}

\author*[1]{\fnm{Theresa} \sur{Wagner}}\email{theresa.wagner@math.tu-chemnitz.de}

\author[2]{\fnm{John W.} \sur{Pearson}}\email{j.pearson@ed.ac.uk}

\author[1]{\fnm{Martin} \sur{Stoll}}\email{martin.stoll@math.tu-chemnitz.de}

\affil[1]{\orgdiv{Department of Mathematics}, \orgname{University of Technology Chemnitz}, \orgaddress{\city{Chemnitz}, \postcode{09107}, \country{Germany}}}

\affil[2]{\orgdiv{School of Mathematics}, \orgname{The University of Edinburgh}, \orgaddress{\city{Edinburgh}, \postcode{EH9 3FD}, \country{UK}}}

\abstract{In this paper we consider the numerical solution to the soft-margin support vector machine optimization problem. This problem is typically solved using the SMO algorithm, given the high computational complexity of traditional optimization algorithms when dealing with large-scale kernel matrices. In this work, we propose employing an NFFT-accelerated matrix--vector product using an ANOVA decomposition for the feature space that is used within an interior point method for the overall optimization problem. As this method requires the solution of a linear system of saddle point form we suggest a preconditioning approach that is based on low-rank approximations of the kernel matrix together with a Krylov subspace solver. We compare the accuracy of the ANOVA-based kernel with the default LIBSVM implementation. We investigate the performance of the different preconditioners as well as the accuracy of the ANOVA kernel on several large-scale datasets.}




\maketitle

\section{Introduction and motivation}
The training of support vector machines (SVMs) leads to large-scale quadratic programs (QPs) \cite{scholkopf2002learning}. An efficient way to solve these optimization problems is the sequential minimal optimization (SMO) algorithm introduced in \cite{platt1998sequential}. The main motivation for the design of the SMO algorithm comes from the fact that existing optimization methods, i.e., quadratic programming approaches, cannot handle the large-scale dense kernel matrix efficiently. The SMO algorithm is motivated by the result obtained in \cite{osuna1997improved} that showed the optimization problem can be decomposed into the solution of smaller subproblems, which avoids the large-scale dense matrices.

When tackling the training as a QP programming task, the use of interior point methods (IPM) has also been studied in the seminal paper of Fine and Scheinberg in \cite{fine2001efficient}: the authors use a low-rank approximation of the kernel matrix and propose a pivoted low-rank Cholesky decomposition to approximate the kernel matrix. A similar matrix approximation was also proposed in \cite{harbrecht2012low}. Fine and Scheinberg then use either a Sherman--Morrison--Woodbury formulation or a Cholesky-based approach for handling the linear algebra within the interior point method in an efficient manner.

The approximation of kernel matrices or the approximation of matrix--vector products with them is an active area of research. One can use fast Gauss transform methods \cite{yang2003improved} or non-equispaced fast Fourier transform methods (NFFT) for approximating the matrix--vector products \cite{alfke2018nfft,nestler2022learning,stoll2020literature}. In this paper, we revisit the application of an interior point method (IPM) \cite{FGW02,Gondzio25,NN94,potra2000interior,Wright} to this problem, and we combine this with fast approaches for the matrix representation of the kernel matrix. We especially focus on the cases when the dimensionality of the feature space $\R^d$ grows for larger $d$, situations where fast matrix--vector product approaches typically do not scale perfectly for an increasing number of data points.

In this paper, we propose the use of a feature grouping approach that creates an ANOVA kernel, only coupling up to three features in an individual kernel summand. The resulting matrix--vector product can be efficiently realized using the NFFT approach. In order to mitigate the slow convergence of iterative solvers for this problem we then apply preconditioning approaches based on low-rank approximations, based on developments in preconditioning such as those in \cite{avron2017faster,cai2022fast,cutajar2016preconditioning,martinsson2016randomized}. These can then be efficiently embedded into our suggested kernel IPM.

Before discussing our approach in more detail we briefly recall the optimization problem at the heart of training support vector machines. Following this, in Section~\ref{sec:ipm} we introduce the interior point method to be employed. In Section~\ref{sec:NFFT} we outline the NFFT technology for matrix--vector products, and in Section~\ref{sec:lr} we discuss approaches for preconditioning the kernel matrix. We then describe our overall preconditioned iterative solver in Section~\ref{sec:Precond}, whereupon we showcase its numerical performance in Section~\ref{sec:NumericalResults} and present concluding remarks in Section~\ref{sec:conclusions}.

\subsection*{The SVM optimization problem}
We assume that we are given $n$ data points $x_i\in\R^d$ associated with a response $y_i\in\left\lbrace -1,1\right\rbrace.$ The goal is to learn the classification boundary separating the two classes encoded in the response $y_i$ via SVMs \cite{cortes1995support}.  

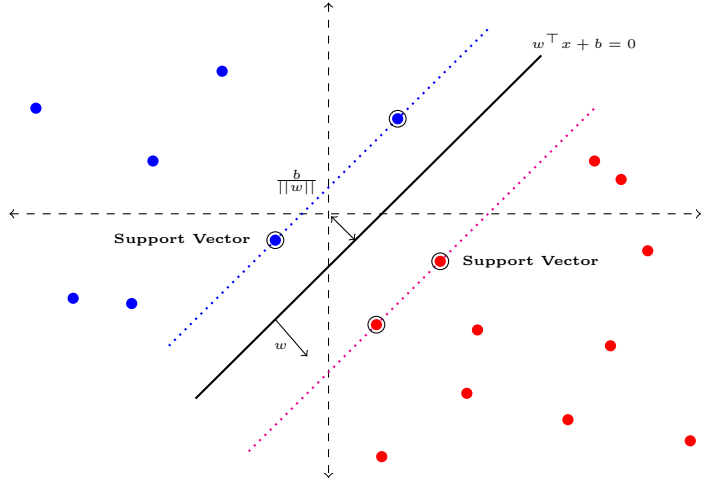
\begin{figure}
\begin{center}
\begin{tikzpicture}[scale=0.70]
\draw[<->, dashed] (-6,0) -- (7,0);
\draw[<->, dashed] (0,4) -- (0,-5);
\draw[blue,thick,dotted] (-3,-2.5) -- (3,3.5);
\draw[thick] (-2.5,-3.5) -- (4,3);
\draw[magenta,thick,dotted] (-1.5,-4.5) -- (5,2);
\draw[<->] (0.05,-0.05) -- (0.5,-0.5);
\draw[->] (-1,-2) -- (-0.4,-2.7);
\fill[red] (0.9,-2.1) circle (3pt);
\draw (0.9,-2.1) circle (4.5pt);
\fill[red] (2.1,-0.9) circle (3pt);
\draw (2.1,-0.9) circle (4.5pt);
\fill[red] (5,1) circle (3pt);
\fill[red] (5.5,0.65) circle (3pt);
\fill[red] (6,-0.7) circle (3pt);
\fill[red] (2.8,-2.2) circle (3pt);
\fill[red] (5.3,-2.5) circle (3pt);
\fill[red] (2.6,-3.4) circle (3pt);
\fill[red] (4.5,-3.9) circle (3pt);
\fill[red] (6.8,-4.3) circle (3pt);
\fill[red] (1,-4.6) circle (3pt);
\fill[blue] (1.3,1.8) circle (3pt);
\draw (1.3,1.8) circle (4.5pt);
\fill[blue] (-1,-0.5) circle (3pt);
\draw (-1,-0.5) circle (4.5pt);
\fill[blue] (-5.5,2) circle (3pt);
\fill[blue] (-3.3,1) circle (3pt);
\fill[blue] (-2,2.7) circle (3pt);
\fill[blue] (-3.7,-1.7) circle (3pt);
\fill[blue] (-4.8,-1.6) circle (3pt);
\node at (-2.75,-0.5) {\textbf{\tiny{Support Vector}}};
\node at (3.8,-0.9) {\textbf{\tiny{Support Vector}}};
\node at (-0.9,-2.5) {\tiny{${w}$}};
\node at (4.8,3.25) {\tiny{${w}^{\top}{x}+b=0$}};
\node at (-0.6,0.6) {\tiny{$\frac{b}{||{w}||}$}};
\end{tikzpicture}
\end{center}
\caption{Illustration for the support vector machine maximal margin classifier.}
\label{fig:svmillu}
\end{figure}

The original method was introduced based on a hard maximal margin separation as illustrated in Figure~\ref{fig:svmillu}.
SVMs allow us to obtain a decision boundary for the classification of the response variables $y_i$.  For binary classification of the given data points $(x_i,y_i)$ with $y_i\in\left\lbrace-1,1\right\rbrace$, a hard constraint is considered in the optimization as
$$
y_i(w^{\top}x_i+b)\geq 1.
$$
As this constraint is designed to not allow for misclassified points, we soften it using the constraint
$$
y_i(w^{\top}x_i+b)\geq 1-\xi_i,\quad \xi_i\geq 0.
$$
For this soft margin constraint, we consider the optimization problem
\begin{align*}
 \min_{w,\xi}\ &\frac{1}{2}w^{\top}w+C\sum_i \xi_i\\
 \nonumber&y_i(w^{\top}x_i+b)\geq 1-\xi_i,\quad \xi_i\geq 0,
\end{align*}
where the parameter $C$ controls the amount of misclassification. With the help of Lagrangian duality and Lagrange multipliers $\alpha_i$ this problem is then transformed to obtain 
\begin{align*}
\max_{\substack{0\leq\alpha_i\leq C}}\sum_{i}\alpha_i-\frac{1}{2}\sum_{i,j}\alpha_iy_i\alpha_jy_jx_i^{\top}x_j,
\end{align*}
a quadratic programming problem with the box constraint $0\leq\alpha_i\leq C$ and the constraint $y^{\top}\alpha=0$. This problem can be formulated in matrix--vector form as
\begin{align*}
\max_{\substack{0\leq\alpha_i\leq C}}e^{\top}\alpha-\frac{1}{2}\alpha^{\top}YXX^{\top}Y\alpha,
\end{align*}
where $e=[1,\ldots,1]^{\top},$ $\alpha=[\alpha_1,\ldots,\alpha_n]$ and $Y=\mathrm{diag}(y_1,\ldots,y_n).$ In particular, we see that the structure of $XX^{\top}$ as
\[
XX^{\top}=
\begin{bmatrix}
x_1^{\top}x_1&\ldots&x_1^{\top}x_n\\
\vdots&\ddots&\vdots\\
x_n^{\top}x_1&\ldots&x_n^{\top}x_n\\
\end{bmatrix}\quad\textrm{ using }~ 
X=
\begin{bmatrix}
- x_1^{\top}-\\
\vdots\\
-x_n^{\top}-\\
\end{bmatrix},\
X^{\top}=
\begin{bmatrix}
x_1&\ldots&x_n\\
\end{bmatrix}
\] 
shows that this is a Gram matrix. The decision boundary for this problem is linear, which is often not sufficient to obtain accurate classification. As such, we consider a formulation suited for data that are not linearly separable. The above derivation can be carried out entirely in a reproducing kernel Hilbert space (RKHS) \cite{scholkopf2002learning} setting, where all inner products are evaluated via a kernel function $\kappa(x_i,x_j).$  The objective function of the optimization then becomes,
\begin{align}
\label{kernelsvm}
\max_{\substack{\alpha}}\sum_{i}\alpha_i-\frac{1}{2}\sum_{i,j}\alpha_iy_i\alpha_jy_j\kappa(x_i,x_j)
= \max_{\substack{\alpha}}e^{\top}\alpha-\frac{1}{2}\alpha^{\top}YKY\alpha,
\end{align}
subject to the constraints $0\leq\alpha_i\leq C$ and $y^{\top}\alpha=0$, with $K_{ij} = \kappa (x_i, x_j )$.  The kernel matrix $K$ is a symmetric matrix that for positive kernels is positive semi-definite. In this paper we consider the case of $K\in \R^{n,n}$ being large-scale, as then the treatment with QP approaches becomes prohibitively expensive. The default approach for solving this optimization problem is the sequential minimal optimization method, which uses the following idea. For problem \eqref{kernelsvm}, the method chooses two Lagrange multipliers $\alpha_i$ and $\alpha_j$ for which the objective function is then maximized while all other $\alpha$ values are kept fixed. This process is repeated until convergence. We refer to \cite{platt1998sequential} for further details. 

Our goal in this paper is to investigate the use of an interior point scheme for solving the kernel SVM problem, by taking advantage of efficient numerical linear algebra methods. We therefore introduce the underlying numerical scheme in Section~\ref{sec:ipm}.
\section{The interior point method}
\label{sec:ipm}
We apply an interior point method to solve the SVM optimization problem. Among many excellent references, we refer to \cite{FGW02,Gondzio25,NN94,potra2000interior,Wright} for outlines of the interior point technology, and in particular to \cite{Gondzio25} for a discussion of the form of the method which we apply here. In such a method, an initial guess is taken which satisfies the bound constraints $0\leq\alpha_i\leq C$, whereupon a barrier sub-problem is solved at every interior point iteration, involving a \emph{barrier parameter} $\mu>0$, where the stationary point of the following Lagrangian is sought:
\begin{equation*}
\ L_{\mu}(\alpha,\lambda)=e^{\top} \alpha-\frac{1}{2}\alpha^{\top} YK Y\alpha+\lambda y^{\top} \alpha+\mu\sum_j \text{log}(\alpha_j)+\mu\sum_j \text{log}(C-\alpha_j),
\end{equation*}
with $Y$ a diagonal matrix containing the entries of $y$, and $e\in\R^n$ the vector of ones. The parameter $\mu$ is progressively reduced towards zero, either by a constant \emph{barrier reduction parameter} $\sigma$, or by a parameter $\sigma_{\texttt{it}}$ which changes at every interior point iteration (denoted $\texttt{it}$).

Upon differentiating $L_{\mu}$ with respect to each entry of $\alpha$ as well as $\lambda$, the following optimality conditions are obtained:
\begin{align*}
\ YK Y\alpha-\lambda y -\mu e./\alpha+\mu e./(Ce-\alpha)={}&e, \\
\ -y^{\top} \alpha={}&0.
\end{align*}

Taking $\bar{\alpha}$ and $\bar{\lambda}$ to be a current interior point iterate for $\alpha$ and $\lambda$, with $\Delta\alpha$ and $\Delta\lambda$ the updates, gives the following Newton conditions for an updated interior point step:
\begin{align*}
\ YK Y\Delta\alpha+{}&\mu\left(e./\bar{\alpha}^2+e./(Ce-\bar{\alpha})^2\right)\circ\Delta\alpha-\Delta\lambda \, y \\
\ &{}=e-YK Y\bar{\alpha}+\bar{\lambda}y +\mu e./\bar{\alpha}-\mu e./(Ce-\bar{\alpha}), \\
\ -y^{\top} \Delta\alpha&{}=y^{\top} \bar{\alpha},
\end{align*}
where $\circ$ corresponds to the component-wise product of two vectors, $\bar{\alpha}^2$ and $(Ce-\bar{\alpha})^2$ denote vectors with the entries squared component-wise, and we apply the MATLAB notation `./' to denote entry-wise division by a vector.

Now, denoting $\Theta$ as a diagonal matrix with the $j$-th entry equal to $\mu(1/\bar{\alpha}_j^2+1/(C-\bar{\alpha}_j)^2)$, the linear system which must be solved to determine the direction of the Newton step is as follows:
\begin{equation}\label{Newton}
\ \underbrace{\left[\begin{array}{cc}
YK Y+\Theta & -y \\ -y^{\top} & 0 \\
\end{array}\right]}_{\mathcal{A}}\left[\begin{array}{c}
\Delta\alpha \\ \Delta\lambda \\
\end{array}\right]=\left[\begin{array}{c}
e-YK Y\bar{\alpha}+\bar{\lambda}y+\mu e./\bar{\alpha}-\mu e./(Ce-\bar{\alpha}) \\ y^{\top}\bar{\alpha} \\
\end{array}\right].
\end{equation}
We address the solution of the system~\eqref{Newton} in Section~\ref{sec:Precond}. Having (inexactly) solved this system, we then determine step-lengths such that the Newton update continues to satisfy the bound constraints. In Algorithm~\ref{alg:IPM} we present the structure of the interior point method which we use to solve the support vector machine problem numerically.


\begin{algorithm}[ht]
  \caption{Interior point method for support vector machines \label{alg:IPM}}
  \begin{algorithmic}[1]
  \State \textbf{Specify Parameters and Initialize Interior Point Method}
  \State $\gamma_0=0.99995$, step-size factor to boundary.
  \State \texttt{tol}, stopping tolerance.
  \State Barrier parameter $\mu_{0}= 1.0$.
  \State Initial guesses for $\alpha^0$, $\lambda^0$, assumed to be such that $\xi_{\alpha}^{0}:=y^{\top} \alpha^{0} \neq 0${, $\xi_{\lambda}^{0}:=e-YK Y\alpha^{0} +\lambda^{0} y +\mu e./\alpha^{0}-\mu e./(Ce-\alpha^{0}) \neq 0$}.
  \State Set iteration number $\texttt{it}=0$.
  \State \textbf{Interior Point Method Loop}
  \While{($\mu>\texttt{tol}$ or $\left\|\xi_{\alpha}^{\texttt{it}}\right\|/\left\|\xi_{\alpha}^0\right\|>\texttt{tol}$ or $\left\|\xi_{\lambda}^{\texttt{it}}\right\|/\left\|\xi_{\lambda}^0\right\|>\texttt{tol}$)}
    \State Reduce barrier parameter $\mu_{\texttt{it}+1}=\sigma_{\texttt{it}}\mu_{\texttt{it}}$.
    \State Solve Newton system \eqref{Newton} for Newton direction $\Delta\alpha$, $\Delta\lambda$.
    \State Find $s_{\alpha}$, $s_{\lambda}$ such that bound constraints on primal and dual variables hold.
    \State Set $s_{\alpha}=\gamma_{0}s_{\alpha}$, $s_{\lambda}=\gamma_{0}s_{\lambda}$.
    \State Make step: $\alpha^{\texttt{it}+1}=\alpha^{\texttt{it}}+s_{\alpha}\Delta\alpha$, $\lambda^{\texttt{it}+1}=\lambda^{\texttt{it}}+s_{\lambda}\Delta\lambda$.
    \State Update infeasibilities: $\xi_{\alpha}^{\texttt{it}+1}=y^{\top} \alpha^{\texttt{it}+1}${, \\ \qquad $\xi_{\lambda}^{\texttt{it}+1}=e-YK Y\alpha^{\texttt{it}+1} +\lambda^{\texttt{it}+1} y +\mu e./\alpha^{\texttt{it}+1}-\mu e./(Ce-\alpha^{\texttt{it}+1})$.}
    \State Set iteration number $\texttt{it}=\texttt{it}+1$.
  \EndWhile
  \end{algorithmic}
\end{algorithm}

\section{NFFT-based matrix--vector products and high dimensional feature spaces}\label{sec:NFFT}
As can be seen from the previous section, the main cost in the interior point approach lies in the solution of the linear system of saddle point form \cite{benzi2005numerical,elman2014finite}. Krylov subspace methods such as MINRES \cite{minres} or GMRES \cite{gmres} rely on matrix--vector multiplications with the matrix $\mathcal{A}$ to build up the Krylov subspace $\mathcal{K}_l(\mathcal{A},r)=\mathrm{span}\left\lbrace r,\mathcal{A}r, \mathcal{A}^2r,\ldots, \mathcal{A}^{l-1}r\right\rbrace$ in terms of the residual vector $r$.
The matrix $\mathcal{A}$, and in particular the kernel matrix $K$, are dense matrices in our framework. In the context of the graph Laplacian $L=D-W$ (cf.\ \cite{stoll2020literature,von2007tutorial}), $W$ is an adjacency matrix, which resembles the structure of the kernel matrix $K$, and in this paper we use the notation $D$ to denote a diagonal matrix, which in this setting is defined as the degree matrix. While in the graph case one could overcome the density of $W$ by sparsifying the graph via $\varepsilon$-neighborhood or $k$-nearest neighborhood approaches, this is not applicable for the kernel case.  For this case, we examine the structure of the matrix--vector product $Kv$ with $v\in\R^n$ and $K$ being the Gaussian kernel function:
$$
\kappa(x_i,x_j)=\exp\left(-\frac{\norm{x_j-x_i}^2}{\ell^2}\right),
$$
where $\ell$ is a scaling parameter for the shape of the Gaussian. With this kernel function, we then obtain the kernel matrix $K$ and want to approximate the matrix--vector product $Kv$, where $v$ is a vector of appropriate dimensionality. An entry of this vector then has the form
\[
\left(Kv\right)_j
=\sum_{i=1}^{n}v_i\exp{\left(-\frac{\norm{x_j-x_i}^2}{\ell^2}\right)}\quad \forall j=1,\ldots,n.
\]
Following \cite{alfke2018nfft} we write this again using the kernel function 
\[
\sum_{i=1}^{n}v_i\exp{\left(-\frac{\norm{x_j-x_i}^2}{\ell^2}\right)}=\sum_{i=1}^{n}v_i\kappa(\tilde{x})
\]
with $\tilde{x}=x_j-x_i.$ The goal is now to approximate the kernel function $\kappa(\cdot)$ using a trigonometric polynomial
$$\tilde{\kappa}(\tilde{x}):=\sum_{\mathbf{J}\in I_N}\tilde{b}_{\mathbb{J}}e^{2\pi i\mathbf{J}\tilde{x}},\quad I_N:=\left\lbrace-\frac{N}{2},-\frac{N}{2}+1,\ldots,\frac{N}{2}-1\right\rbrace,$$
where $N\in 2\mathbb{N}$ is the bandwidth and $\tilde{b}_{\mathbb{J}}$ are the Fourier coefficients. We can then rewrite the component of the matrix--vector product as
\begin{align*}
\left(Kv\right)_j
&=\sum_{i=1}^{n}v_i\kappa(\tilde{x})\approx\sum_{i=1}^{n}v_i\tilde{\kappa}(\tilde{x})\\
&=\sum_{i=1}^{n}v_i\sum_{\mathbf{J}\in I_N}\tilde{b}_{\mathbb{J}}e^{2\pi i\mathbf{J}\tilde{x}}\\
&=\sum_{i=1}^{n}v_i\sum_{\mathbf{J}\in I_N}\tilde{b}_{\mathbb{J}}e^{2\pi i\mathbf{J}(x_j-x_i)}\\
&=\sum_{\mathbf{J}\in I_N}\tilde{b}_{\mathbb{J}}\left(\sum_{i=1}^{n}v_ie^{-2\pi i\mathbf{J}x_i}\right)e^{2\pi i\mathbf{J}x_j}.
\end{align*}
As shown in \cite{alfke2018nfft}, we can rely on the efficient evaluation of both the inner and the outer sum using the NFFT \cite{potts2003fast} method, here such that the cost of one matrix--vector multiplication becomes 
$\mathcal{O}(m^dn+N^d\log{N})$, where $m$ is a fixed window parameter. As observed in \cite{alfke2018nfft} this provides us with a method of linear complexity for fixed accuracy, which without further assumptions is typically limited to moderate dimensions where $d\leq 3.$
For the applications that we are aiming for in this paper we cannot limit ourselves to this case, and we wish to adjust the NFFT approach for the case $d>3.$

In a similar spirit, \cite{morariu2008automatic} uses the improved fast Gauss transform for reducing the cost of the matrix--vector multiplication. Previous comparisons of both approaches \cite{alfke2018nfft} revealed that the performances are in a comparable range. For higher dimensional feature spaces $d>10$, the authors in \cite{march2015askit} suggest an approach based on skeleton approximations. While this approach is not employed in this paper we believe this could provide interesting further topics for future investigations.

To use the NFFT for speeding up the solution of the linear system within an interior point method, a special kernel function is required that fulfills several requirements. First, the chosen kernel must be able to approximate the kernel function well using a trigonometric polynomial, which is not straightforward. Next, the input dimension for the kernel and hence for every NFFT computation should not exceed three, so that the computational advantage of the NFFT can be exploited. Our aim is to construct a kernel that is based on a sum of kernels that each only rely on at most three features.  At the same time, as many feature interactions as possible shall be involved and the number of kernels should still be kept low. For this we use the extended Gaussian ANOVA kernel
\begin{align*}
      \kappa ( x_i, x_j) = \sum_{l=1}^P \eta_l \underbrace{\exp \left( - \frac{\| x_j^{\mathcal W_l} - x_i^{\mathcal W_l} \|^2}{\ell_l^2} \right)}_{\kappa_l \left( x_i^{\mathcal{W}_l}, x_j^{\mathcal{W}_l} \right)}.
\end{align*}
Let us now explain the ingredients of this kernel. The number of kernels $P$ satisfies {$P \leq \lceil \frac{d}{3} \rceil$}. For selection of the relevant features to be included in each kernel we use the index sets $\mathcal W_l = \{ w_1^l, w_2^l, w_3^l\} \subseteq \{1, \dots, d\}$, which are determined following the feature's mutual information score ranking. Based on the \emph{windows of features} $\mathcal{W}_l$, we obtain the data points restricted to those features as
\begin{align*}
    x_i^{\mathcal W_l} = \begin{bmatrix} x_i^{w_1^l} & x_i^{w_2^l} & x_i^{w_3^l} \end{bmatrix}^\top
    \quad \text{and} \quad
    x_j^{\mathcal W_l} = \begin{bmatrix} x_j^{w_1^l} & x_j^{w_2^l} & x_j^{w_3^l} \end{bmatrix}^\top,
\end{align*}
where {$i,j = 1, \dots, n$} and $l = 1, \dots, P$. Note that each kernel $\kappa_l$ is assigned an individual kernel parameter $\ell_l$. For more details on this kernel we refer to \cite{nestler2022learning}. We are now able to combine multiple kernels, each relying on at most three features, so that the NFFT can be applied to approximate the action of the kernel matrices corresponding to an overall matrix--vector product $(K_1+\ldots+K_P)v=K_1v+\ldots+K_Pv$, where $K_{l_{ij}} = \kappa_l (x_i^{\mathcal{W}_l}, x_j^{\mathcal{W}_l})$. We embed this into the interior point method and the corresponding iterative solver within this. We next need to discuss the preconditioning approach that is required to ensure fast convergence of the iterative solver, as unpreconditioned solvers for the system
$$
\mathcal{A}x=g
$$
often converge very slowly. We therefore consider the (left or right) preconditioned system
$$
\mathcal{P}^{-1}\mathcal{A}x=\mathcal{P}^{-1}g \quad \text{or} \quad \mathcal{A}\mathcal{P}^{-1}(\mathcal{P}x)=g,
$$
where the choice of the preconditioner $\mathcal{P}$ is discussed in Section~\ref{sec:Precond}.

\section{Preconditioning the kernel matrix\label{sec:lr}}
We wish to investigate the use of low-rank approaches for preconditioning the kernel matrix $K$ or its approximations based on the ANOVA NFFT. We briefly recall several approaches that have been suggested for approximating $K$. Note that we here first discuss the approximation of a single kernel matrix and later comment on the sum of several kernel matrices.
\subsection*{Pivoted Cholesky decomposition}
The Cholesky decomposition of $K=LL^{\top}$ is a staple of scientific computing, computational statistics, and optimization. Unfortunately, the decomposition comes at cubic cost given the dimensionality of the kernel matrix. As a remedy, one can use a pivoted Cholesky decomposition of small rank as a low-rank approximation to the kernel matrix. This idea was suggested by Fine and Scheinberg in \cite{fine2001efficient}. The pivoted Cholesky decomposition is discussed in some detail in \cite{GoluVanl96} and we follow the derivation there. It was also proposed in \cite{harbrecht2012low} for applications to different kernel matrices. We briefly recall the method here and point out that it is derived from the outer product form of the full Cholesky decomposition. For this form we decompose the matrix as follows:
$$
K=
\begin{bmatrix}
K_{11}&K_{12}\\
K_{21} & K_{22}
\end{bmatrix}=
\begin{bmatrix}
\sqrt{K_{11}}&0\\
K_{21}/\sqrt{K_{11}}&I_{n-1}
\end{bmatrix}
\begin{bmatrix}
1&0\\
0&K_{22}-K_{21}K_{12}/K_{11}
\end{bmatrix}
\begin{bmatrix}
\sqrt{K_{11}}&K_{12}/\sqrt{K_{11}}\\
0&I_{n-1}
\end{bmatrix}
$$
with $K_{11} \in \R^+$, $K_{21} \in \R^{n-1}$, and $K_{22} \in \R^{n-1,n-1}$. As the matrix $K$ is positive (semi-)definite, the square root of $K_{11}$ is well defined and one can show that the matrix $K_{22}-K_{21}K_{12}/K_{11}$ also has the same property. Using this relationship one can inductively show that the Cholesky decomposition exists and can be computed in the style of a Gaussian elimination (cf.\ \cite[Alg. 4.2.2]{GoluVanl96}).
For rank-deficient matrices, or cases when we are interested in a low-rank approximation of the kernel matrix $K$, we proceed with a pivoting strategy based on the symmetric transformation $\Pi K \Pi^{\top},
$ with $\Pi$ being a permutation matrix. Note that such a permutation will keep diagonal elements on the diagonal. We give the detailed method in Algorithm~\ref{alg:pchol} to obtain the approximation
$$
K\approx LL^{\top}.
$$
\begin{algorithm}[ht]
  \caption{Pivoted Cholesky decomposition of rank $k$ (cf.\ \cite[Alg. 4.2.4]{GoluVanl96}) \label{alg:pchol}}
  \begin{algorithmic}[1]
  \State$r=0$
  \For{$j=1,\ldots,n$}
    \State Find $l \geq j$ such that $l=\argmax(K_{j,j},\ldots,K_{n,n})$.
    \If{$K_{l,l}>0$ and $r<k$}
    \State $r=r+1$
    \State Swap $K_{:,l}$ and $K_{:,j}$, as well as $K_{l,:}$ and $K_{j,:}.$
    \State Set $K_{j,j}=\sqrt{K_{j,j}}.$
    \State Set $K_{j+1:n,j}=K_{j+1:n,j}/K_{j,j}.$
    \For{$i=j+1,\ldots,n$}
    \State $K_{i:n,i}=K_{i:n,i}-K_{i:n,j}K_{i,j}$
    \EndFor
    \EndIf
  \EndFor
  \State Define $L$ as the lower-triangular part of $K$ up to column $k$.
  \end{algorithmic}
\end{algorithm}

Recently, the authors of \cite{chen2022randomly} have introduced a variant of the pivoted Cholesky decomposition coined \textit{randomly-pivoted Cholesky decomposition.} One can use a greedy selection process, where the element with the largest absolute value on the diagonal of the updated matrix is selected as the pivot element. One can also select the pivoting elements at random, and the authors in \cite{chen2022randomly} select each pivot as being sampled in proportion to the diagonal entries of the current residual matrix.


\subsection*{Nystr\"om approximation}
The Nystr\"om method is a well-known technique in machine learning for the approximation of the kernel matrix \cite{drineas2005nystrom}. We briefly introduce the method here, as based on \cite{martinsson2019randomized}. The basic idea is the following approximation:
$$
K\approx(KQ)(Q^{\top}KQ)^{-1}(KQ)^{\top},
$$
where the matrix $Q\in\R^{n,k}$ has orthogonal columns. The solution of the linear system resulting from the $(Q^{\top}KQ)^{-1}$ term can be carried out by employing a Cholesky decomposition, but we found this to be unstable at times. Instead, once this matrix has eigenvalues close to zero, we use a $LDL^{\top}$ decomposition where we set the diagonal entries of $D$ to a fixed value as soon as they are smaller than a predefined threshold. If $Q$ is constructed from $k$ columns of a permutation matrix, one obtains the original Nystr\"om method, where the kernel matrix  
$$
K=
\begin{bmatrix}
K_{11}&K_{12}\\
K_{21}&K_{22}
\end{bmatrix}
$$
is structured. The block $K_{11}\in\R^{k,k}$ contains the $k$ components of $K$ that we use for the approximation, and $K_{12}\in\R^{k,n}$ their interactions with the remaining data points. Note that for simplicity of presentation we have assumed that $K$ is permuted so the relevant components appear in the upper left corner. The goal of the Nystr\"om method is to avoid the explicit storage and computation of $K_{22}\in\R^{n-k,n-k}.$  We thus use 
$$
K\approx
\begin{bmatrix}
K_{11}&K_{12}\\
K_{21}&K_{21}K_{11}^{-1}K_{12}
\end{bmatrix}.
$$
In the more general setup, $Q$ is obtained from the relation $Q = \mathrm{orth}(KG)$, where $G\in\R^{n,k}$ is a Gaussian matrix with normally distributed random entries and
$\mathrm{orth}$ denotes column-wise orthonormalization. Note that we only require matrix--vector products with the matrix $K$ to compute this Nystr\"om approximation, and that given an approximation to the matrix--vector product we do not require $K$ explicitly.
\subsection*{Random Fourier features}
The idea of the random Fourier features approach \cite{rahimi2007random} is to use a set of feature space basis functions such as
$$
\phi(x,w,b)=\cos(w^{\top}x+b),
$$
where we draw different values $(w_i,b_i)$ randomly from some distribution. This is done on the basis of Bochner's theorem \cite{plonka2018numerical}, where for a real-valued kernel function $\kappa(x,x')=\kappa(x-x')$ with $x,x'\in\R^d$ there is a probability distribution $\rho(\cdot)$ over $\R^d$, which we assume to be such that
$$
\kappa(x,x')=\int_{\R^d}\rho(w)\cos(w^{\top}(x-x'))dw.
$$
Given $w$ that is distributed based on $\rho(\cdot)$ we know that 
$$
\kappa(x,x')=\mathbb{E}\left[\cos(w^{\top}(x-x'))\right].
$$
Drawing $k$ samples from the distribution $\rho(\cdot)$, i.e.,
$w_1,\ldots,w_k\sim\rho(w)$, we obtain
$$
\frac{1}{k}\sum_{j=1}^{k}\cos(w_j^{\top}(x-x'))\approx\mathbb{E}\left[\cos(w^{\top}(x-x'))\right]=\kappa(x,x').
$$
We know that $\kappa(x,x')=\left\langle\phi(x),\phi(x')\right\rangle,$ where the function $\phi$ is given by
$\phi_i(x)=\sqrt{2} \, \cos(w_i^{\top}x+b_i)$
and $b_i$ are uniformly sampled from the interval $0$ and $2\pi$ (cf.\ \cite{rahimi2008weighted}). If we now collect all the $\phi_i$ into one vector, we obtain
$$
\phi(x)=\frac{1}{\sqrt{k}}\left[\phi_1(x),\ldots,\phi_k(x)\right]^{\top}\in\R^{k}
$$
and an approximation to the kernel matrix via 
\begin{align*}
K&=
\begin{bmatrix}
k(x_1,x_1)&\ldots&k(x_1,x_n)\\
\vdots&\ddots&\vdots\\
k(x_n,x_1)&\ldots&k(x_n,x_n)\\
\end{bmatrix}\\
&\approx
\begin{bmatrix}
\phi(x_1)^{\top} \phi(x_1)&\ldots&\phi(x_1)^{\top} \phi(x_n)\\
\vdots&\ddots&\vdots\\
\phi(x_n)^{\top} \phi(x_1)&\ldots&\phi(x_n)^{\top} \phi(x_n)\\
\end{bmatrix}\\
&=
\underbrace{
\begin{bmatrix}
\phi(x_1)^{\top}\\
\vdots\\
\phi(x_n)^{\top}\\
\end{bmatrix}}_{\in\R^{n,k}}
\underbrace{
\begin{bmatrix}
\phi(x_1),
\hdots,
\phi(x_n)\\
\end{bmatrix}}_{\in\R^{k,n}},
\end{align*}
which gives the low-rank approximation using random Fourier features.

\subsection*{Some comments}
We have now seen three approaches where the kernel matrix is approximated via a low-rank approximation of the form 
$$
K\approx Z^{\top}Z
$$
for the purpose of preconditioning. Note that the above approaches require only the knowledge of very few entries of the kernel matrix to obtain the low-rank approximation. There are several methods that can be embedded into the IPM solver via providing a black-box matrix--vector product. Here we use the non-equispaced fast Fourier transform method \cite{alfke2018nfft,nestler2022learning} but one could also utilize the fast Gauss transform~\cite{yang2003improved, elgammal2003efficient}. If in the Nyström method $Q$ does not represent columns of the identity matrix, we require the fast evaluation of $KQ$ via the NFFT scheme. It would also be possible to compute other low-rank approximations such as the randomized singular value decomposition \cite{halko2011finding}.

Once we want to apply the above low-rank techniques to the ANOVA formulation 
$$
K_1+\ldots+K_P,
$$
we assume that the low-rank approaches are applied to each individual matrix via 
$$
K_1+\ldots+K_P\approx Z_1^{\top}Z_1+\ldots+Z_P^{\top}Z_P 
$$
which we can write as
$$
\underbrace{
\begin{bmatrix}
Z_1^{\top}&Z_2^{\top}&\ldots&Z_P^{\top}
\end{bmatrix}}_{\in\R^{n,Pk}}
\underbrace{\begin{bmatrix}
Z_1\\
Z_2\\
\vdots\\
Z_P
\end{bmatrix}}_{\in\R^{Pk,n}}.
$$
This is again a low-rank approximation to the sum of matrices, but with an increased dimensionality. We next discuss how to handle the efficient evaluation of the overall preconditioner.

\section{Preconditioned iterative solver}\label{sec:Precond}
It is clear that the main computational burden of Algorithm~\ref{alg:IPM} is the solution of the linear system~\eqref{Newton}. As such, solving this system efficiently is a key consideration in this paper. We focus on iterative methods based on Krylov subspaces \cite{minres,saad2003iterative,gmres}. These methods are the most efficient when combined with a suitable preconditioner and we here rely on well-established saddle point theory \cite{benzi2005numerical,MGW}, based on which we wish to apply a preconditioner $\mathcal{P}$ for $\mathcal{A}$ of the form:
\begin{equation*}
\ \mathcal{P}=\left[\begin{array}{cc}
\widehat{A} & 0 \\ -y^{\top} & -1 \\
\end{array}\right].
\end{equation*}
Here $\widehat{A}$ denotes an approximation to the matrix $A:=YKY+\Theta$, that is the $(1,1)$-block of $\mathcal{A}$. When $A$ is approximated exactly, we may readily verify that
\begin{equation*}
\ \mathcal{P}^{-1}\mathcal{A}=\left[\begin{array}{cc}
I & -A^{-1}y \\ 0 & y^{\top} A^{-1}y \\
\end{array}\right],
\end{equation*}
and hence the preconditioned system has $n$ eigenvalues equal to $1$, with the remaining eigenvalue equal to $y^{\top} A^{-1}y$. This shows that the preconditioner $\mathcal{P}$ is a good candidate for approximating $\mathcal{A}$, provided an appropriate method for approximating $A$ may be devised. Due to the non-symmetry of $\mathcal{P}$, we apply this within the GMRES algorithm \cite{gmres}, which also allows for a non-symmetric approximation of $A$. It is of course possible to apply an analogous symmetric positive definite preconditioner within the MINRES algorithm \cite{minres}, for which one may guarantee a specific convergence rate based on eigenvalues of the preconditioned system, as long as $\widehat{A}$ is itself symmetric positive definite.

As discussed in Section~\ref{sec:lr}, the kernel matrix is approximated by the ANOVA NFFT and the preconditioner will be used for this approximation instead. To apply the saddle point preconditioner, it is clear that the key step is to apply the action of $\widehat{A}^{-1}$ to a generic vector. This involves using $A$, a matrix which consists of the sum of the diagonal barrier matrix $\Theta$, which is positive definite and typically becomes more ill-conditioned as the interior point method progresses, and the matrix $YKY$.




Having applied one of the low-rank approximations we obtain
\begin{equation*}
 \widehat{A}=
 \Theta+YZ^{\top}ZY
 =
 Y\left[Y^{-1}\Theta Y^{-1}+Z^{\top}Z\right]Y.
\end{equation*}
Note that $D:=Y^{-1}\Theta Y^{-1}$ is a readily-computable diagonal matrix at each interior point iteration.

Now, using the Sherman--Morrison--Woodbury identity \cite{Woodbury}, we may write
\begin{align*}
\ \widehat{A}^{-1}={}&Y^{-1}\left[D+Z^{\top}Z\right]^{-1}Y^{-1} \\
={}&Y^{-1}\left[D^{-1}-D^{-1}Z^{\top}(I_k+ZD^{-1}Z^{\top})^{-1}ZD^{-1}\right]Y^{-1}\\
={}&Y^{-1}D^{-1}\left[I-Z^{\top}(I_k+ZD^{-1}Z^{\top})^{-1}ZD^{-1}\right]Y^{-1}.
\end{align*}
Substituting in $D^{-1}:=Y\Theta^{-1} Y$,
\begin{align*}
Y^{-1}Y\Theta^{-1} Y\left[I-Z^{\top}(I_k+ZY\Theta^{-1} YZ^{\top})^{-1}ZY\Theta^{-1} Y\right]Y^{-1}\\
=\Theta^{-1} -\Theta^{-1}YZ^{\top}(I_k+ZY\Theta^{-1} YZ^{\top})^{-1}ZY\Theta^{-1}.
\end{align*}
The main computational workload which needs to be undertaken \emph{at each interior point iteration} is to compute the $k\times k$ matrix $I_k+ZD^{-1}Z^{\top},$ which we will do using a direct solver such as an LU decomposition readily available in Python. The overall algorithm for applying the preconditioner
$$
\mathcal{P}\left[\begin{array}{c}
   x_1 \\ x_2 \\
\end{array}\right]=\left[\begin{array}{c}
  g_1 \\ g_2 \\
  \end{array}\right]
$$
requires first the solution with $\widehat{A}$, which in turn needs the LU decomposition of the $k\times k$ matrix $I_k+ZD^{-1}Z^{\top}$ once per IPM iteration, and all other operations with $D$ and $Y$ are computationally trivial as these matrices are diagonal. Once $x_1$ is computed we can update the last entry $x_2=-y^{\top} x_1-g_2.$ We have thus obtained an efficient preconditioning scheme, which we now test on several challenging examples.

\section{Numerical results}\label{sec:NumericalResults}

In this section we present the numerical results for the proposed method on benchmark datasets. The corresponding implementations are available in the GitHub repository NFFTSVMipm, see \url{https://github.com/wagnertheresa/NFFTSVMipm}. All experiments were run on a computer with $8 \, \times$ Intel Core $\text{i}7-7700$ CPU @ $3.60$ GHz processors with NV106 graphics and $16.0$ GiB of RAM.

Before starting the training routine it is important to preprocess the data to prevent unwanted effects. Our implementations cover balancing the training data and $z$-score normalization. By doing so we prevent the algorithm from simply predicting the over-represented class and ensure all data is standardized, which is necessary for applying the NFFT approach. Note that the statistics for the scaling are solely computed with respect to the training data to prevent train--test-contamination. The train--test-split is $0.5$.

We consider the UCI datasets HIGGS ($n = 11000000$, $d = 28$)~\cite{higgs_dataset} and SUSY ($n = 5000000$, $d = 18$)~\cite{susy_dataset}, as well as the LIBSVM dataset cod-rna ($n = 488565$, $d = 8$)~\cite{codrna_dataset}.

In the following subsections we demonstrate the performance of our NFFTSVMipm method with several preconditioners by analyzing their setup time and the number of GMRES iterations required per IPM step. Moreover, we compare its runtime and predictive power to the results of LIBSVM~\cite{chang2011libsvm}, a state-of-the-art library for solving SVMs. With this we provide an extensive understanding of NFFTSVMipm's overall performance for different experimental settings. All results presented in this section were generated with the parameter choices listed in Table~\ref{tab:params} unless stated otherwise. Deciding when to stop the IPM iteration was linked to the performance of the preconditioned GMRES method within the IPM scheme.

\begin{table}[h]
\centering
\begin{tabular}{|r|r l|}
 \hline
\multicolumn{3}{| c |}{General Parameter Setting}\\
 \hline
 Maximum number of interior point iterations & $\texttt{iter}_{\text{ip}}=$ & \hspace{-1.1em}$50$ \\
 IPM convergence tolerance & $\texttt{tol}_{\text{ip}}=$ & \hspace{-1.1em}$10^{-1}$ \\
 Maximum number of GMRES iterations & $\texttt{iter}_{\text{GMRES}}=$ & \hspace{-1.1em}$100$ \\
 GMRES convergence tolerance & $\texttt{tol}_{\text{GMRES}}=$ & \hspace{-1.1em}$10^{-3}$ \\
 IPM barrier reduction parameter & $\sigma=$ & \hspace{-1.1em}$0.6$ \\
 Initial step size within IPM & $\texttt{gap} =$ & \hspace{-1.1em}$0.99995$ \\
 FastAdjacency setup & $\texttt{fastadj}_{\text{setup}}=$ & \hspace{-1.1em}``default'' \\ \hline
 \multicolumn{3}{| c |}{Parameter Setting for Preconditioners} \\
 \hline
 Maximum error tolerance for Cholesky precond. & $\texttt{err}_{\text{tol}}=$ & \hspace{-1.1em}$10^{-5}$ \\ \hline
 \multicolumn{3}{| c |}{Parameter Bounds for Random Search} \\
 \hline
 Kernel parameter $\ell$ & $\ell \in $ & \hspace{-1.1em}$[ 10^{-1}, 10 ]$ \\
 Relative weight of error vs. margin, such that $0 \leq \alpha \leq C$ & $C \in$ & \hspace{-1.1em}$[0.1, 0.7]$ \\ \hline
\end{tabular}
\caption{Parameter setting for the experiments presented in this paper.}
\label{tab:params}
\end{table}

\subsection{Comparison of low-rank preconditioners}

First we want to examine the effect of realizing NFFTSVMipm with different preconditioners of various respective ranks. For this we implement the preconditioners introduced in Section~\ref{sec:lr}, within the framework described in Section~\ref{sec:Precond}. In the following we abbreviate the greedy-based and randomized pivoted Cholesky preconditioners with Cholesky (greedy) and Cholesky (rp), respectively, and the random Fourier features preconditioner with RFF. We denote the number of training data points and the number of test data points by $n_{\text{train}}$ and $n_{\text{test}}$, respectively. 

\begin{figure}
    \centering
    \includegraphics[width=\textwidth]{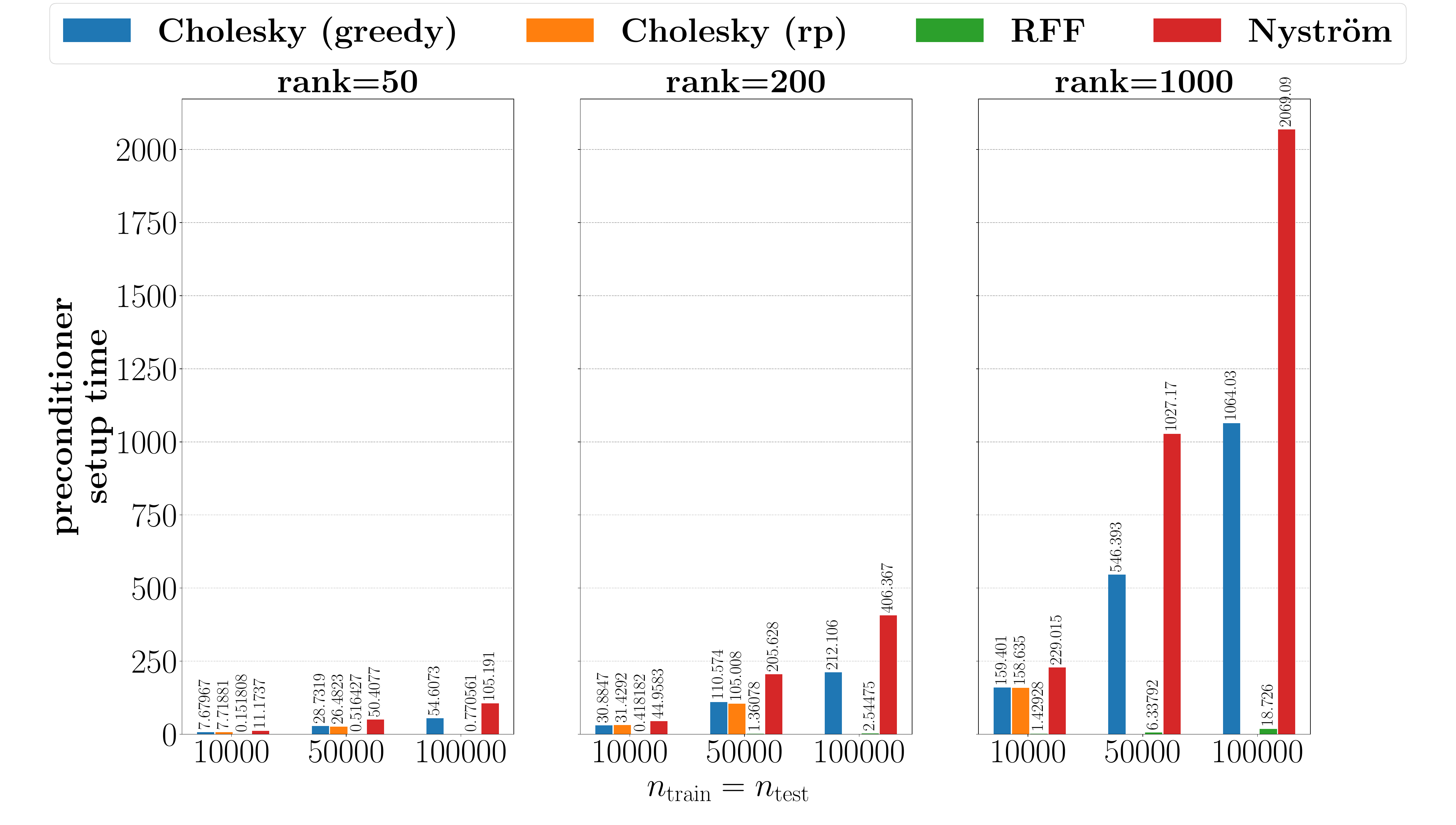}
    \caption{Preconditioner setup time in seconds for different preconditioners of various ranks for the SUSY dataset performed with the NFFTSVMipm.}
    \label{fig:susy_precond_timing}
\end{figure}

Figure~\ref{fig:susy_precond_timing} visualizes the setup time for these preconditioners for various ranks and subset sizes. Throughout the entire experiment both Cholesky preconditioners take roughly the same setup time. By comparison, the calculation of the Nyström preconditioner takes considerably longer, which is all the more significant for higher preconditioner ranks. The green bars representing the setup time for the RFF preconditioner can barely be seen, which is due to its computation consisting solely of generating random samples and evaluating trigonometric expressions. This leads to incredibly low computation times. Note that some results for the Cholesky (rp) preconditioner are missing for larger subsets and a large preconditioner rank, since the sum of the diagonal entries of the decomposition matrix reaches zero in those cases. The algorithm then aborts, because it involves dividing by this value. Therefore, this method is very sensitive to the ratio between small subset sizes $n_{\text{train}}$ and the rank of the preconditioner. As expected, the higher the values of the rank and the subset size, the larger is the setup time.

We are looking for a suitable preconditioner to embed into our approach, targeted at outperforming state-of-the-art methods in terms of their computational complexity. Even though the preconditioner setup time is of relevance it is decisive that the preconditioner enables fast convergence. Thus, we examine the mean number of GMRES iterations required per IPM step next. This serves as a good indicator of how well the preconditioned method solves the Newton system within an IPM and therewith of the overall prediction quality.

\begin{figure}
    \centering
    \includegraphics[width=\textwidth]{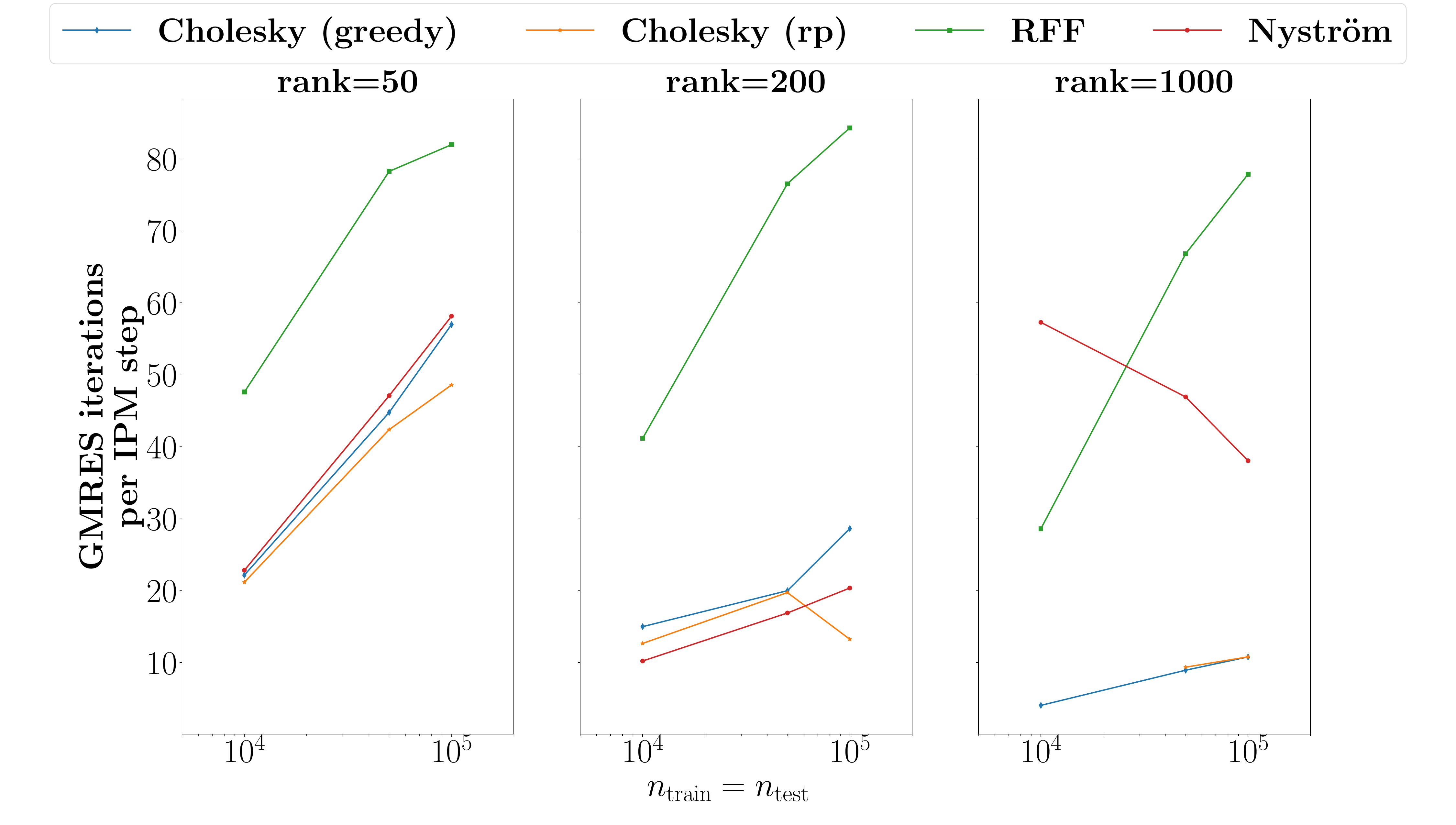}
    \caption{Mean number of GMRES iterations per IPM step for different preconditioners of various ranks for the SUSY dataset performed with the NFFTSVMipm.}
    \label{fig:susy_precond_rank}
\end{figure}

As can been seen in Figure~\ref{fig:susy_precond_rank}, the mean number of GMRES iterations within the IPM method required with both Cholesky and the Nyström preconditioners are roughly within the same range throughout the experiment. Note that the Cholesky (rp) preconditioner aborts for rank $1000$ again for the same reason as described earlier. Naturally, the higher the preconditioner rank the smaller the number of GMRES iterations required within the IPM steps. This is due to the preconditioners of higher rank yielding a better approximation of the dense kernel matrices and hence leading to $\mathcal{P}$ serving as a better preconditioner for $\mathcal{A}$ in system~\eqref{Newton}. Merely for rank $1000$, the Nyström preconditioner seems to improve for growing $n_{\text{train}}$. By contrast, the RFF preconditioner greatly exceeds the iteration numbers of the competitive preconditioners. In most cases it requires nearly twice as many iterations per IPM step, which leads to a significantly worse convergence behavior and ultimately a longer training time for the SVM. Moreover, depending on the exact problem and parameter setting, the IPM is more likely to not converge leading to a worse prediction quality. 

Having taken into account the performance regarding several aspects of the  different preconditioners, we realize that the RFF preconditioner is not well suited for our problem. Having shown promising results regarding the setup time it however cannot keep up with the precipitated convergence behavior yielded by the competitive preconditioners. The Nyström preconditioner, on the contrary, cannot compete in terms of the setup time and does not show a stable convergence behavior. This leaves us with the Cholesky preconditioners where both show similar behavior regarding setup time and number of GMRES iterations. However, the greedy method results in increased robustness. It can be seen that the Cholesky (rp) preconditioner fails for a range of values in Figure~\ref{fig:susy_precond_timing} and for a different range of values in Figure~\ref{fig:susy_precond_rank}. We believe that an adaptive strategy might resolve this issue but for our purposes we choose the Cholesky (greedy) preconditioner for the remainder of these results.

\subsection{Comparison of convergence tolerances}

\begin{figure}
    \centering
    \includegraphics[width=\textwidth]{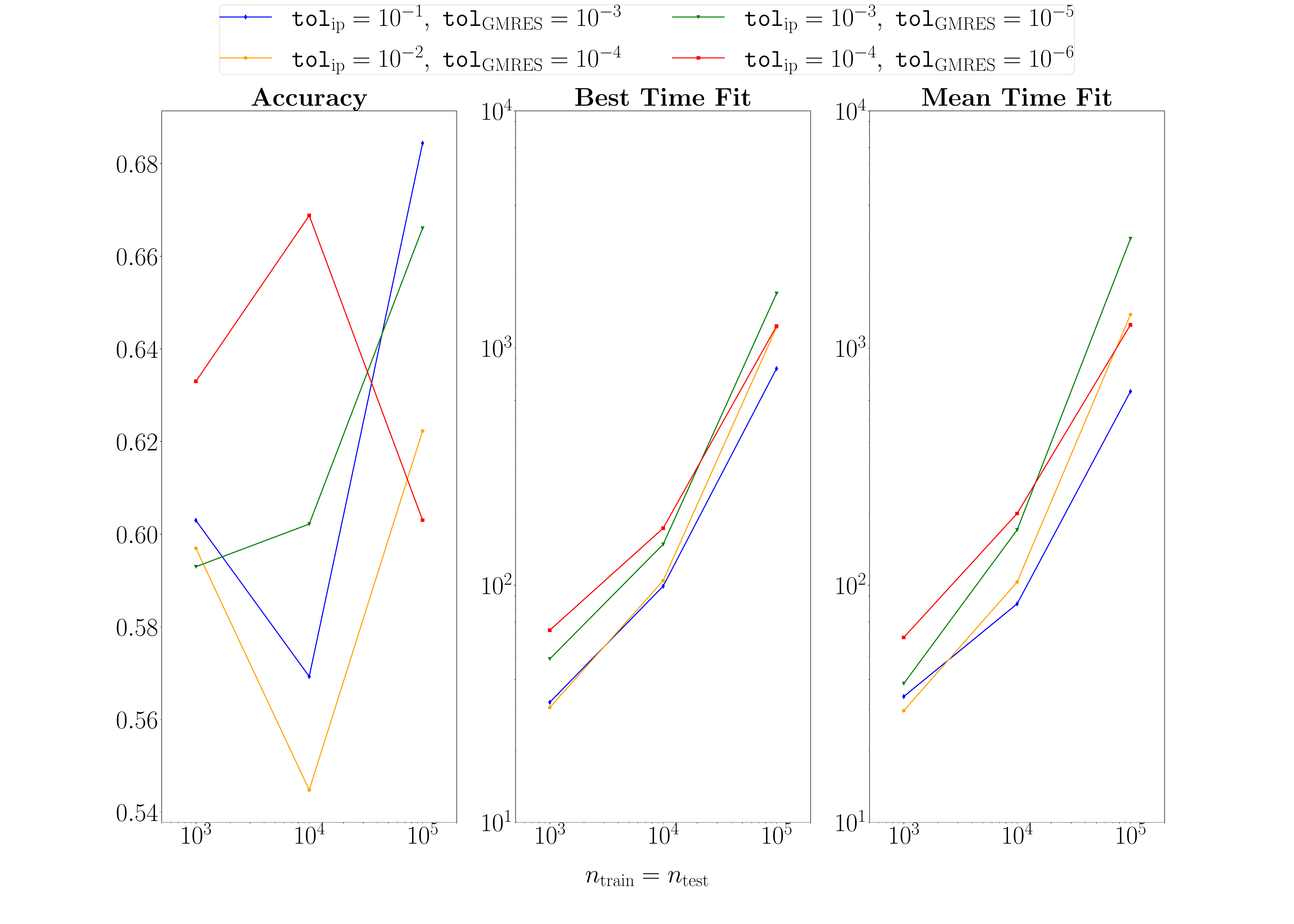}
    \caption{Performance of SVM trained with NFFTSVMipm (with Cholesky (greedy) preconditioner of rank $200$) for different IPM and GMRES convergence tolerances for the HIGGS dataset.}
    \label{fig:tols}
\end{figure}

In Figure~\ref{fig:tols} we illustrate for the HIGGS dataset how our method performs on different convergence tolerances for both IPM and GMRES. It can be seen that a tighter tolerance for the IPM results in an improved accuracy for the largest values of $n_{\text{train}}$ and $n_{\text{test}}$. As expected the computing time for fitting the parameters increases, but not dramatically so. 

\subsection{Comparison with LIBSVM}

\begin{figure}[ht!]
    \centering
    \includegraphics[width=\textwidth]{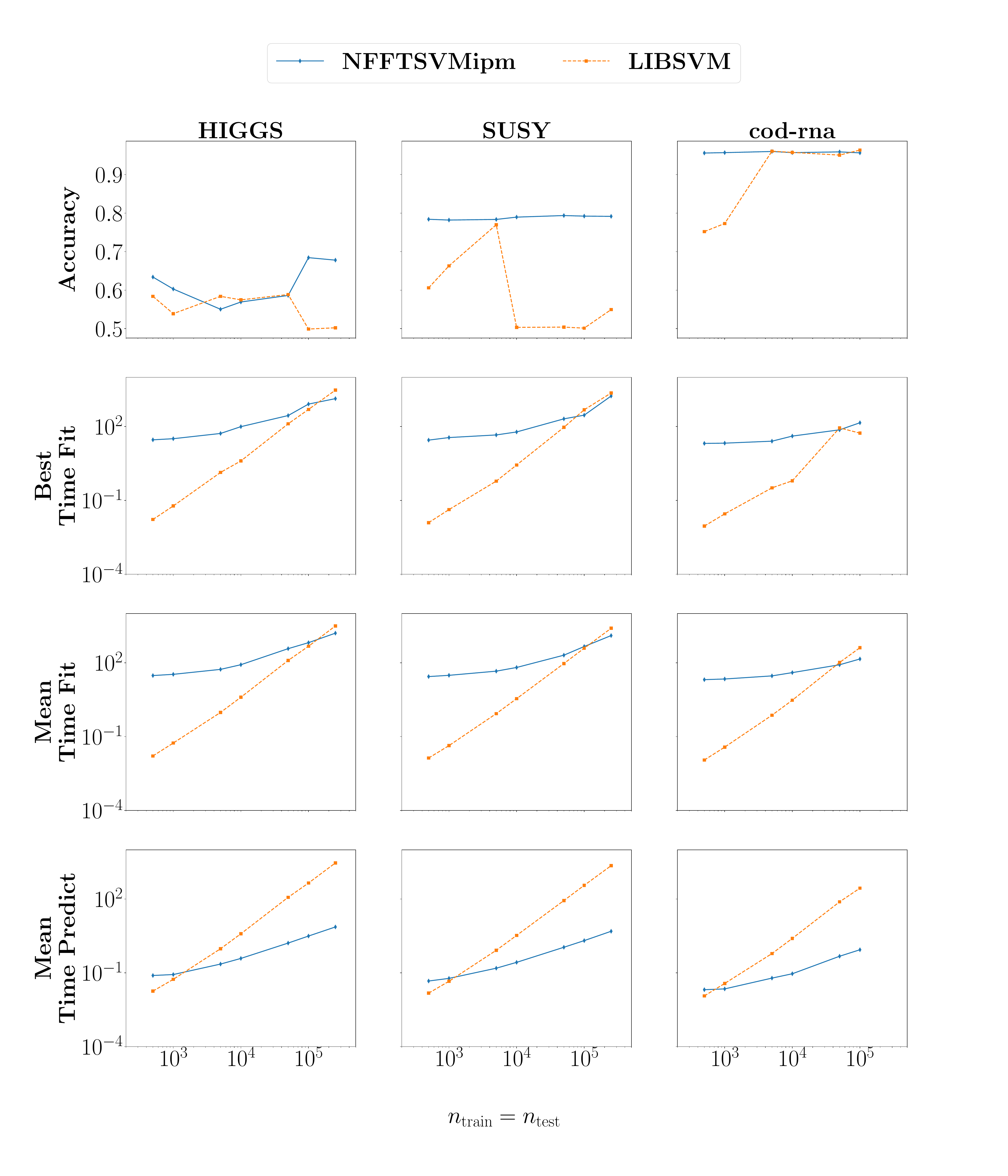}
    \caption{Performance of SVM trained with LIBSVM and our approach NFFTSVMipm (with Cholesky (greedy) preconditioner of rank $200$) on benchmark datasets.}
    \label{fig:results}
\end{figure}

Having analyzed the suitability of various preconditioners for our problem in the previous subsection, we now compare the proposed approach in terms of runtime and accuracy with the Python implementation of LIBSVM \cite{chang2011libsvm}, which relies on the SMO algorithm. In order to examine both techniques, we use a random search strategy for selecting the most successful parameter combinations. We run both methods for $25$ randomly chosen parameter combinations within previously defined bounds, see Table~\ref{tab:params}. The ANOVA setup, due to the inclusion of several kernel functions with individual kernel parameters $\ell$ for the corresponding windows, has a much larger parameter space than the standard SVM. We report the accuracy, the runtime to obtain the best-performing parameters, as well as the mean runtime for fitting the model and for prediction.

As can be seen from Figure~\ref{fig:results}, the computational advantage of training the SVM with the NFFTSVMipm method is evident for large datasets for more than $100000$ training data points. This is observed consistently for all benchmark datasets considered in this experiment with the feature dimension ranging between $8$ and $28$, i.e., with a number of feature windows $P$ between $3$ and $10$. It can be seen that our NFFTSVMipm method already reaches the break-even point at $50000$ training data points, with the cod-rna dataset. In fact this dataset has the smallest feature dimension, leading to fewer feature windows, and thus to fewer NFFT approximations within every IPM iteration. By contrast, when comparing the computing time for unseen test data the NFFTSVMipm overtakes the LIBSVM routine already for much smaller subset sizes. The prediction times are similar for $1000$ test data points, and when working with $100000$ data points the NFFTSVMipm method outperforms the LIBSVM by several orders of magnitude. This is due to the prediction routine relying on far fewer kernel--vector multiplications and therefore NFFT evaluations than the training process. 

The last three rows in Figure~\ref{fig:results} show the runtimes depicted as log--log plots to better see the difference in slopes between the two approaches. It stands out that NFFTSVMipm's curves start at a higher level than LIBSVM's. Comparing the slope for both methods indicates that NFFTSVMipm eventually outperforms LIBSVM and then provides a computational advantage. While our proposed method cannot keep up with LIBSVM's runtimes for small subset sizes, since it requires some NFFT setup time, NFFTSVMipm is generally targeted at big data and shows promising results regarding accuracy. 

Depending on the data set and exact problem setting, NFFTSVMipm yields a up to $25$ percent better accuracy than LIBSVM. While LIBSVM takes into account all feature information at once within one kernel matrix, NFFTSVMipm works with an additive kernel approach that splits up the feature dimensions into smaller windows of features and therefore covers lower-order feature interactions.

The potential advantage in predictive performance when working with an additive kernel approach is certainly a pressing research question. However, this exceeds the scope of this paper and will be a subject of future research.

\section{Conclusions}
\label{sec:conclusions}

In this paper we have analyzed the applicability of NFFT-accelerated kernel--vector products to SVM problems. Moreover we designed a preconditioner for the saddle point system consisting of a low-rank approximation of the kernel matrix and a Krylov subspace solver. We presented several options of such low-rank preconditioners for the kernel and examined their performance. Finally, we compared the proposed NFFTSVMipm method to the state-of-the-art LIBSVM solver. We illustrated that in the regime of large-scale datasets the NFFTSVMipm approach provides a computational advantage, highlighting the advantages of making use of the linear algebra structure of the problem, which leads to us not requiring access to the full (storage-intensive) kernel matrix. In particular, we observed that separating the features into multiple kernels via the ANOVA decomposition can provide possibly better accuracy, in some cases significantly so. It remains to investigate this effect in future work. We emphasize that utilizing the ANOVA decomposition and fast matrix--vector product could also exhibit a benefit when solving SVMs using other methods, for instance the alternating direction method of multipliers (ADMM) \cite{GM75,GM76} or other first-order methods, which would also be a fruitful avenue for further investigation.

\section*{Acknowledgments}
TW and MS gratefully acknowledge their support from the Bundesministerium f\"{u}r Bildung und Forschung (BMBF) grant 01\,$\vline$\,S20053A (project SA$\ell$E). JWP gratefully acknowledges financial support from the Engineering and Physical Sciences Research Council (EPSRC) grant EP/S027785/1. 


\end{document}